\magnification=1200
\def\qed{\unskip\kern 6pt\penalty 500\raise -2pt\hbox
{\vrule\vbox to 10pt{\hrule width 4pt\vfill\hrule}\vrule}}
\centerline{DIFFERENTIATING THE ABSOLUTELY CONTINUOUS}
\centerline{INVARIANT MEASURE OF AN INTERVAL MAP $f$}
\centerline{WITH RESPECT TO $f$.}
\bigskip
\centerline{by David Ruelle\footnote{*}{Mathematics Dept., Rutgers University, and IHES.  91440 Bures sur Yvette, France.\break $<$ruelle@ihes.fr$>$}.}
\bigskip\bigskip\noindent
	{\leftskip=2cm\rightskip=2cm\sl {\bf Abstract.}  Let the map $f:[-1,1]\to[-1,1]$ have a.c.i.m. $\rho$ (absolutely continuous $f$-invariant measure with respect to Lebesgue).  Let $\delta\rho$ be the change of $\rho$ corresponding to a perturbation $X=\delta f\circ f^{-1}$ of $f$.  Formally we have, for differentiable $A$, 
$$   \delta\rho(A)=\sum_{n=0}^\infty\int\rho(dx)\,X(x){d\over dx}A(f^nx)   $$
but this expression does not converge in general.  For $f$ real-analytic and Markovian in the sense of covering $(-1,1)$ $m$ times, and assuming an {\it analytic expanding} condition, we show that 
$$	\lambda\mapsto\Psi(\lambda)=\sum_{n=0}^\infty\lambda^n
	\int\rho(dx)\,X(x){d\over dx}A(f^nx)      $$
is meromorphic in ${\bf C}$, and has no pole at $\lambda=1$.  We can thus formally write $\delta\rho(A)=\Psi(1)$.
\par}
\vfill\eject
	We postpone a discussion of the significance of our result, and start to describe the conditions under which we prove it.  Note that these conditions are certainly too strong: suitable differentiability should replace analyticity, and a weaker Markov property should be sufficient.  But the point of the present note is to show how it is that $\Psi(\lambda)$ has no pole at $\lambda=1$, rather than deriving a very general theorem.  
\medskip\noindent
{\bf Setup.}
\medskip
	We assume that $f:[-1,1]\to[-1,1]$ is real analytic and piecewise monotone on $[-1,1]$ in the following sense: there are points $c_j$ ($j=0,\ldots,m$, with $m\ge2$) such that $-1=c_0<c_1<\ldots<c_{m-1}<c_m=1$ and, for $j=0,\ldots,m$, 
$$	f(c_j)=(-1)^{j+1}      $$
We assume that on $[-1,1]$ the derivative $f'$ vanishes only on $Z=\{c_1,\ldots,c_{m-1}\}$, and that $f''$ does not vanish on $Z$.  For $j=1,\ldots,m$, we have $f[c_{j-1},c_j]=[-1,1]$.  In particular, $f$ is Markovian.  We shall also assume that $f$ is {\it analytically expanding} in the sense of Assumption A below.  The purpose of this note is to prove the following:
\medskip\noindent
{\bf Theorem.}  {\sl Under the above conditions, and Assumption A stated later, there is a unique $f$-invariant probability measure $\rho$ absolutely continuous with respect to Lebesgue on $[-1,1]$.  If $X$ is real-analytic on $[-1,1]$, and $A\in{\cal C}^1[-1,1]$, then 
$$	\Psi(\lambda)
=\sum_{n=0}^\infty\lambda^n\int_{-1}^1\rho(dx)\,X(x){d\over dx}A(f^nx)      $$
extends to a meromorphic function in ${\bf C}$, without pole at $\lambda=1$.}
\medskip
	Our proof depends on a change of variable which we now explain.  We choose a holomorphic function $\omega$ from a small open neighborhood $U_0$ of $[-1,1]$ in ${\bf C}$ to a small open neighborhood $W$ of $[-1,1]$ in a Riemann surface which is 2-sheeted over ${\bf C}$ near $-1$ and $1$.  We call $\varpi=\omega^{-1}:W\to U_0$ the inverse of $\omega$.  We assume that $\omega(-x)=-\omega(x)$, $\omega(\pm1)=\pm1$, $\omega[-1,1]=[-1,1]$, $\omega'(\pm1)=\omega'''(\pm1)=0$.  We have thus 
$$	\omega(\pm(1-\xi))=\pm(1-C\xi^2+D\xi^4\ldots)      $$
with $C>0$ and, if $a>0$, 
$$	\varpi(\pm(1-a\xi^2+b\xi^3\ldots))
	=\pm(1-\sqrt{{a\over C}}\xi+{b\over2\sqrt{aC}}\xi^2\ldots)      $$
[We may for instance take 
$$	\omega(x)=\sin{\pi x\over2}\qquad,
	\qquad\varpi(x)={2\over\pi}\arcsin x      $$
or 
$$	\qquad\omega(x)={1\over16}(25x-10x^3+x^5)\qquad,
	\qquad\varpi(x)={16\over25}x\ldots\qquad\qquad]      $$
\indent
	The function $g:\varpi\circ f\circ\omega$ from $[-1,1]$ to $[-1,1]$ has monotone restrictions to the intervals $\varpi[c_{j-1},c_j]=[d_{j-1},d_j]$.  It is readily seen that $g_j$ extends to a holomorphic function in a neighborhood of $[d_{j-1},d_j]$, and that 
$$	g_1(-1+\xi)=-1+\sqrt{f'(-1)}\xi+\alpha_-\xi^3\ldots      $$
$$	g_m(1-\xi)=(-1)^{m+1}(1-\sqrt{|f'(1)|}\xi-\alpha_+\xi^3\ldots)      $$
with no $\xi^2$ terms in the right-hand sides [this follows from our choice of $\omega$, which has no $\xi^3$ term].  One also finds that, for $j=1,\ldots,m-1$ 
$$	g_j(d_j-\xi)
=(-1)^{j+1}(1-\sqrt{|f''(c_j)|\over2C}\omega'(d_j)\xi+\gamma_j\xi^2\ldots)   $$
$$	g_{j+1}(d_j+\xi)
=(-1)^{j+1}(1-\sqrt{|f''(c_j)|\over2C}\omega'(d_j)\xi-\gamma_j\xi^2\ldots)   $$
where $\gamma_j$ is the same in the two relations.  We note the following easy consequences of the above developments:
\medskip\noindent
{\bf Lemma 1.}  {\sl Let $\psi_j:[-1,1]\to[d_{j-1},d_j]$ be the inverse of $g_j$ for $j=1,\ldots,m$ (increasing for $j$ odd, decreasing for $j$ even).  Then 
$$	\psi_1(-1+\xi)=-1+{1\over\sqrt{f'(-1)}}\xi+\beta_-\xi^3      $$
$$    \psi_m((-1)^{m+1}(1-\xi))=1-{1\over\sqrt{|f'(1)|}}\xi+\beta_+\xi^3    $$
(there are no $\xi^2$ terms in the right-hand sides).  If $j<m$, 
$$	\psi_j((-1)^{j+1}(1-\xi))
=d_j-\sqrt{2C\over|f''(c_j)|}{1\over\omega'(d_j)}\xi+\delta_j\xi^2      $$
$$	\psi_{j+1}((-1)^{j+1}(1-\xi))
=d_j+\sqrt{2C\over|f''(c_j)|}{1\over\omega'(d_j)}\xi+\delta_j\xi^2      $$
(with the same coefficient $\delta_j$).}\qed  
\medskip
	As inverses of the $g_j$, the functions $\psi_j$ extend to holomorphic functions on a neighborhood of $[-1,1]$.  We impose now the condition that $f$ is {\it analytically expanding} in the following sense: 
\medskip\noindent
{\bf Assumption A}  {\sl We have $[-1,1]\subset U\subset{\bf C}$, with $U$ bounded open connected, such that the $\psi_j$ extend to continuous functions $\bar U\mapsto{\bf C}$, holomorphic in $U$, and with $\psi_j\bar U\subset U$.  [$\bar U$ denotes the closure of $U$].}
\medskip
	Let $\phi$ be holomorphic on a neighborhood of $\bar U$.  Given a sequence ${\bf j}=(j_1,\ldots,j_\ell,\ldots)$ we define $\phi_{{\bf j}\ell}=\phi\circ\psi_{j_1}\cdots\circ\psi_{j_\ell}$ and note that the $\phi_{{\bf j}\ell}$ are uniformly bounded in a neighborhood of $\bar U$.  We may thus choose $\ell(r)$ for $r=1,2\ldots$ such that the subsequence $(\phi_{{\bf j}\ell(r)})_{r=1}^\infty$ converges uniformly on $\bar U$ to a limit $\tilde\phi_{\bf j}$.  Writing $\tilde U=\cup_{j=1}^m\psi_j\bar U$ we have 
$$	\max_{z\in\bar U}|\phi_{{\bf j}\ell(r)}|
	\ge\max_{z\in\tilde U}|\phi_{{\bf j}\ell(r)}|
	\ge\max_{z\in\bar U}|\phi_{{\bf j}\ell(r+1)}|      $$
so that $\max_{z\in\bar U}|\tilde\phi_{\bf j}|=\max_{z\in\tilde U}|\tilde\phi_{\bf j}|$ and, since $\tilde U$ is compact $\subset U$ connected, $\tilde\phi_{\bf j}$ is constant.  Therefore $\phi$ is constant on $\cap_{\ell=0}^\infty\psi_{j_1}\circ\cdots\circ\psi_{j_\ell}\bar U$.  Since this is true for all $\phi$, the intersection $\cap_{\ell=0}^\infty\psi_{j_1}\circ\cdots\circ\psi_{j_\ell}\bar U$ consists of a single point $\tilde z({\bf j})$.  Given $\epsilon>0$ we can thus, for each ${\bf j}$, find $\ell$ such that ${\rm diam}\psi_{j_1}\circ\cdots\circ\psi_{j_\ell}\bar U<\epsilon$.  Hence (using the compactness of the Cantor set of sequences ${\bf j}$) one can choose $L$ so that the $m^L$ sets 
$$	\psi_{j_1}\circ\cdots\psi_{j_L}\bar U      $$
have diameter $<\epsilon$.  The open connected set 
$$	V=\cup_{j_1,\ldots,j_L}\psi_{j_1}\circ\cdots\psi_{j_L}U      $$
satisfies $[-1,1]\subset V\subset U$, and $\psi_j\bar V=\cup_{j_1,\ldots,j_L}\psi_j\circ\psi_{j_1}\circ\cdots\circ\psi_{j_L}\bar U\subset\cup_{j_0,j_1,\ldots,j_{L-1}}\psi_{j_0}\circ\psi_{j_1}\circ\psi_{i_{L-1}}U=V$.  This shows that $U$ can be replaced in Assumption A by a set $V$ contained in an $\epsilon$-neighborhood of $[-1,1]$.  
\medskip
	Since we have shown above that ${\rm diam}\psi_{j_1}\circ\cdots\psi_{j_L}\bar U<\epsilon$, we see that $\psi_1^L$ maps a small circle around $-1$ strictly inside itself.  We have thus $\psi'_1(-1)<1$ ({\it i.e.}, $f'(-1)>1$) and similarly, if $m$ is odd, $\psi'_m(1)<1$ ({\it i.e.}, $f'(1)>1$).
\medskip
	The following two lemmas state some easy facts to be used later.
\medskip\noindent
{\bf Lemma 2.}  {\sl Let $H$ be the Hilbert space of functions $\bar U\to{\bf C}$ which are square integrable (with respect to Lebesgue) and holomorphic in $U$.  The operator ${\cal L}$ on $H$ defined by 
$$  ({\cal L}\Phi)(z)=\sum_{j=1}^m(-1)^{j+1}\psi'_j(z)\Phi(\psi_j(z))  $$
is holomorphy improving.  In particular ${\cal L}$ is compact and trace-class.}\qed
\medskip\noindent
{\bf Lemma 3.}  {\sl On $[-1,1]$ we have 
$$	({\cal L}\Phi)(x)=\sum_j|\psi'_j(x)|\Phi(\psi_j(x))      $$
hence $\Phi\ge0$ implies ${\cal L}\Phi\ge0$ (${\cal L}$ preserves positivity) and 
$$	\int_{-1}^1dx\,({\cal L}\Phi)(x)=\int_{-1}^1dx\,\Phi(x)      $$
(${\cal L}$ preserves total mass).}\qed
\medskip\noindent
{\bf Lemma 4.}  {\sl ${\cal L}$ has a simple eigenvalue $\mu_0=1$ corresponding to an eigenfunction $\sigma_0>0$.  The other eigenvalues $\mu_k$ ($k\ge1$) satisfy $|\mu_k|<1$, and their (generalized) eigenfunctions $\sigma_k$ satisfy $\int_{-1}^1dx\,\sigma_k(x)=0$.}
\medskip
	Let $(\mu_k,\sigma_k)$ be a listing of the eigenvalues and generalized eigenfunctions of the trace-class operator ${\cal L}$.  For each $\mu_k$ there is some $\sigma_k$ such that ${\cal L}\sigma_k=\mu_k\sigma_k$, hence 
$$	|\mu_k|\int_{-1}^1dx\,|\sigma_k(x)|=\int_{-1}^1dx\,|\mu_k\sigma_k(x)|
	=\int_{-1}^1dx\,|({\cal L}\sigma_k)(x)|      $$
$$  \le\int_{-1}^1dx\,({\cal L}|\sigma_k|)(x)=\int_{-1}^1dx\,|\sigma_k(x)|  $$
hence $|\mu_k|\le1$.  Denote by $S_<$ and $S_1$ the spectral spaces of ${\cal L}$ corresponding to eigenvalues $\mu_k$ with $|\mu_k|<1$, and $|\mu_k|=1$ respectively.  If $\sigma_k\in S_<$ then, for some $n\ge1$, 
$$	0=\int_{-1}^1dx\,(({\cal L}-\mu_k)^n\sigma_k)(x)
	=\int_{-1}^1dx\,(1-\mu_k)^n\sigma_k(x)      $$
hence $\int_{-1}^1dx\,\sigma_k(x)=0$.
\medskip
	On the finite dimensional space $S_1$, there is a basis of eigenvectors $\sigma_k$ diagonalizing ${\cal L}$ (if ${\cal L}|S_1$ had non-diagonal normal form, $||{\cal L}^n|S_1||$ would tend to infinity with $n$, in contradiction with $\int_{-1}^1dx\,|({\cal L}^n\Phi)(x)|\le\int_{-1}^1dx\,|\Phi(x)|$).  We shall now show that, up to multiplication by a constant $\ne0$, we may assume $\sigma_k\ge0$.  If not, because $\sigma_k$ is continuous and the intervals $\psi_{j_1}\circ\cdots\circ\psi_{j_n}[-1,1]$ are small for large $n$ (mixing), we would have $|({\cal L}^n\sigma_k)(x)|<({\cal L}^n|\sigma_k|)(x)$ for some $n$ and $x$.  This would imply $\int_{-1}^1dx\,|({\cal L}^n\sigma_k)(x)|<\int_{-1}^1dx\,|\sigma_k(x)|$ in contradiction with ${\cal L}\sigma_k=\mu_k\sigma_k$ and $|\mu_k|=1$.  From $\sigma_k\ge0$ we get $\mu_k=1$, and the corresponding eigenspace is at most one dimensional (otherwise it would contain functions not $\ge0$).  But we have $1\notin S_<$ because $\int_{-1}^1dx
 \!
,1\ne0$, so that $S_1\ne\{0\}$.  Thus $S_1$ is spanned by an eigenfunction, which we call $\sigma_0$, to the eigenvalue $\mu_0=1$.  Finally, $\sigma_0>0$ because if $\sigma_0(x)=0$ we would have also $\sigma_0(y)=0$ whenever $g^n(y)=x$, which is not compatible with $\sigma_0$ continuous $\ne0$.\qed
\medskip\noindent
{\bf Lemma 5.}  {\sl If we normalize $\sigma_0$ by $\int_{-1}^1dx\,\sigma_0(x)=1$, then $\sigma_0(dx)=\sigma_0(x)dx$ is the unique $g$-invariant probability measure absolutely continuous with respect to Lebesgue on $[-1,1]$.  In particular, $\sigma_0(dx)$ is ergodic.}
\medskip
	For continuous $A$ on $[-1,1]$ we have 
$$	\int_{-1}^1\sigma_0(dx)(A\circ g)(x)=\int_{-1}^1dx\,\sigma_0(x)A(g(x))
=\int_{-1}^1dx\,({\cal L}\sigma_0)(x)A(x)=\int_{-1}^1\sigma_0(dx)A(x)      $$
so that $\sigma_0(dx)$ is $g$-invariant.  Let $\tilde\sigma(x)dx$ be another $g$-invariant probability measure absolutely invariant with respect to Lebesgue.  Then, if $\tilde\sigma\ne\sigma_0$ 
$$	\int_{-1}^1dx\,|\sigma_0(x)-\tilde\sigma(x)|
	=\int_{-1}^1dx\,|({\cal L}(\sigma_0-\tilde\sigma))(x)|      $$
$$	<\int_{-1}^1dx\,({\cal L}|\sigma_0-\tilde\sigma|)(x)
	=\int_{-1}^1dx\,|\sigma_0(x)-\tilde\sigma(x)|      $$
by mixing: contradiction.\qed
\medskip\noindent
{\bf Lemma 6.}  {\sl Let $H_1\subset H$ consist of those functions $\Phi$ with derivatives vanishing at $\pm1$: $\Phi'(-1)=\Phi'(1)=0$.  Then ${\cal L}H_1\subset H_1$ and $\sigma_0\in H_1$}.
\medskip
	${\cal L}H_1\subset H_1$ is an easy calculation using Lemma 1.  Furthermore, by Lemma 4, $\sigma_0=\lim_{n\to\infty}{\cal L}^n{1\over2}$, and ${1\over2}\in H_1$ implies $\sigma_0\in H_1$.\qed
\medskip
	The image $\rho(dx)=\rho(x)dx$ of $\sigma_0(x)dx$ by $\omega$ is the unique $f$-invariant probability measure absolutely continuous with respect to Lebesgue on $[-1,1]$.  We have 
$$	\rho(x)=\sigma_0(\varpi x)\varpi'(x)      $$
Consider now the expression 
$$	\Psi(\lambda)
=\sum_{n=0}^\infty\lambda^n\int_{-1}^1\rho(dx)\,X(x){d\over dx}A(f^nx)      $$
where we assume that $X$ extends to a holomorphic function in a neighborhood of $[-1,1]$ and $A\in{\cal C}^1[-1,1]$.  For sufficiently small $|\lambda|$, the series defining $\Psi(\lambda)$ converges.  Writing $B=A\circ\omega$ and $x=\omega y$ we have 
$$  X(x){d\over dx}A(f^nx)=X(\omega y){1\over\omega'(y)}{d\over dy}B(g^ny)  $$
hence 
$$	\Psi(\lambda)=\sum_{n=0}^\infty\lambda^n\int_{-1}^1dy\,
	\sigma_0(y){X(\omega y)\over\omega'(y)}{d\over dy}B(g^ny)      $$
Defining $Y(y)=\sigma_0(y)X(\omega y)/\omega'(y)$, we see that $Y$ extends to a function holomorphic in a neighborhood of $[-1,1]$, which we may take to be $U$, except for simple poles at $-1$ and $1$.  We may write 
$$   \int_{-1}^1dy\,\sigma_0(y){X(\omega y)\over\omega'(y)}{d\over dy}B(g^ny)
	=\int_{-1}^1dy\,Y(y)g'(y)\cdots g'(g^{n-1}y)B'(g^ny)      $$
$$	=\int_{-1}^1ds\,({\cal L}_0^nY)(s)B'(s)      $$
where 
$$	({\cal L}_0\Phi)(s)=\sum_{j=1}^m(-1)^{j+1}\Phi(\psi_js)      $$
and we have thus 
$$	\Psi(\lambda)=\sum_{n=0}^\infty\lambda^n\int_{-1}^1ds\,
	({\cal L}_0^nY)(s)B'(s)      $$
\noindent
{\bf Lemma 7.}  {\sl Let $H_0\subset H$ be the space of functions vanishing at $-1$ and $1$.  Then ${\cal L}_0H_0\subset H_0$.}
\medskip
	This follows readily from Lemma 1.\qed
\medskip\noindent
{\bf Lemma 8.}  {\sl There are meromorphic functions $\Phi_\pm$ with Laurent series 
$$	\Phi_\pm(z)={1\over z\mp1}+O(z\mp1)      $$
at $\pm1$ and $\Phi_\pm(\mp1)=0$ such that 
$$	{\cal L}_0\Phi_-=\sqrt{f'(-1)}\Phi_-      $$
$$	\left\{
   \matrix{{\cal L}_0\Phi_+=\sqrt{f'(1)}\Phi_+\qquad\hbox{if $m$ is odd}\cr
   {\cal L}_0(\Phi_+/\sqrt{|f'(1)|}+\Phi_-/\sqrt{f'(-1)})=\tilde Y\in H_0
	\qquad\hbox{if $m$ is even}\cr}\right.      $$}
\indent
	Define 
$$	p_\pm(z)={1\over z\mp1}-{1\over4}(z\mp1)      $$
then Lemma 1 yields 
$$	({\cal L}_0-\sqrt{f'(-1)})p_-=u_-\in H_0      $$
$$	\left\{
   \matrix{({\cal L}_0-\sqrt{f'(1)})p_+=u_+\in H_0\qquad\hbox{if $m$ is odd}\cr
   {\cal L}_0p_++\sqrt{|f'(1)|}p_-=u_0\in H_0\qquad\hbox{if $m$ is even}\cr}
	\right.      $$
Since $f'(-1)>1$, Lemma 4 shows that ${\cal L}-\sqrt{f'(-1)}$ is invertible on $H$, hence there is $v_-$ such that 
$$	({\cal L}-\sqrt{f'(-1)})v_-=u'_-      $$
and since $\int_{-1}^1dx\,u'_-(x)=0$, also $\int_{-1}^1dx\,v_-(x)=0$ and we can take $w_-\in H_0$ such that $w'_-=v_-$.  Then 
$$	(({\cal L}_0-\sqrt{f'(-1)})w_-)'=({\cal L}-\sqrt{f'(-1)})w'_-
	=({\cal L}-\sqrt{f'(-1)})v_-=u'_-      $$
so that 
$$	({\cal L}_0-\sqrt{f'(-1)})w_-=u_-      $$
without additive constant because the left-hand side is in $H_0$ by Lemma 7.  In conclusion 
$$	({\cal L}_0-\sqrt{f'(-1)})(p_--w_-)=0      $$
and we may take $\Phi_-=p_--w_-$.  
\medskip
	If $m$ is odd, $\Phi_+$ is handled similarly.  If $m$ is even, taking $\Phi_+=p_+$ and writing $\tilde Y=u_0/\sqrt{|f'(1)|}-w_-$ we obtain
$$	{\cal L}_0({\Phi_+\over\sqrt{|f'(1)|}}+{\Phi_-\over\sqrt{f'(-1)}})
	=\tilde Y\in H_0      $$
which completes the proof.\qed
\medskip
	We have $\sigma_0\in H_1$ (Lemma 6), and $X\circ\omega\in H_1$ by our choice of $\omega$.  Also 
$$	\omega'(\pm(1-\xi))=2C\xi-4D\xi^3\ldots      $$
so that 
$$	Y={\bf C}\Phi_-+{\bf C}\Phi_++H_0      $$
\indent
	If $m$ is odd let $Y=c_-\Phi_-+c_+\Phi_++Y_0$, with $Y_0\in H_0$.  Then 
$$\Psi(\lambda)={c_-\over{1-\lambda\sqrt{f'(-1)}}}\int_{-1}^1ds\,\Phi_-(s)B'(s)
	+{c_+\over{1-\lambda\sqrt{f'(1)}}}\int_{-1}^1ds\,\Phi_+(s)B'(s)
	+\Psi_0(\lambda)      $$
where $\Psi_0$ is obtained from $\Psi$ when $Y$ is replaced by $Y_0$.
\medskip
	If $m$ is even let $Y=c_-\Phi_-+\tilde c({\Phi_+/\sqrt{|f'(1)|}}+{\Phi_-/\sqrt{f'(-1)}})+Y_0$, with $Y_0\in H_0$.  Then 
$$	\Psi(\lambda)
	={c_-\over{1-\lambda\sqrt{f'(-1)}}}\int_{-1}^1ds\,\Phi_-(s)B'(s)
	+\tilde c\int_{-1}^1ds\,({\Phi_+\over\sqrt{|f'(1)|}}
	+{\Phi_-\over\sqrt{f'(-1)}})B'(s)      $$
$$	+\lambda\tilde\Psi(\lambda)+\Psi_0(\lambda)      $$
where $\tilde\Psi(\lambda)$ is obtained from $\Psi$ when $Y$ is replaced by $\tilde Y$.  
\medskip
	Writing $\mu_\pm=\sqrt{f'(\pm1)}$ we see that $\Psi(\lambda)$ has two poles at $\mu_\pm^{-1}$ if $m$ is odd, and one pole at $\mu_-^{-1}$ if $m$ is even; the other poles are those of $\Psi_0(\lambda)$ and possibly $\tilde\Psi(\lambda)$.  Since $Y_0\in H_0$ and ${\cal L}_0H_0\subset H_0$, we have 
$$	\Psi_0(\lambda)
	=\sum_{n=0}^\infty\lambda^n\int_{-1}^1ds\,({\cal L}_0^nY_0)(s)B'(s)
=-\sum_{n=0}^\infty\lambda^n\int_{-1}^1ds\,({\cal L}_0^nY_0)'(s)B(s)      $$
$$   =-\sum_{n=0}^\infty\lambda^n\int_{-1}^1ds\,({\cal L}^nY'_0)(s)B(s)   $$
It follows that $\Psi_0(\lambda)$ extends meromorphically to ${\bf C}$ with poles at the $\mu_k^{-1}$.  We want to show that the residue of the pole at $\mu_0^{-1}=1$ vanishes .  By Lemma 4, $\int_{-1}^1dx\,\sigma_k(x)=0$ for $k\ge1$.  Thus, up to normalization, the coefficient of $\sigma_0$ in the expansion of $Y'_0$ is 
$$	\int_{-1}^1dx\,Y'_0(x)=Y_0(1)-Y_0(-1)=0      $$
because $Y_0\in H_0$.  Therefore $\Psi_0(z)$ is holomorphic at $z=1$, and the same argument applies to $\tilde\Psi(z)$, concluding the proof of the theorem.\qed
\medskip\noindent
{\bf Discussion.}
\medskip
	It can be argued that the {\it physical measure} describing a physical dynamical system is an SRB (Sinai-Ruelle-Bowen) measure $\rho$ (see the recent reviews [11], [2] which contain a number of references), or an a.c.i.m. $\rho$ in the case of a map of the interval.  But, typically, physical systems depend on parameters, and it is desirable to know how $\rho$ depends on the parameters ({\it i.e.}, on the dynamical system).  The dependence is smooth for uniformly hyperbolic dynamical systems (see [5], [6] and references given there), but discontinuous in general.  
\medskip
	The present note is devoted to an example in support of an idea put forward in [8]: that derivatives of $\rho(A)$ with respect to parameters can be meaningfully defined in spite of discontinuities.  An ambitious project would be to have Taylor expansions on a large set $\Sigma$ of parameter values and, using a theorem of Whitney [10], to connect these expansions by a function extrapolating $\rho(A)$ smoothly outside of $\Sigma$.  In a different dynamical situation, that of KAM tori, a smooth extension \`a la Whitney has been achieved by Chierchia and Gallavotti [3], and P\"oschel [4].  
\medskip
	In our study we have considered only a rather special set $\Sigma$ consisting of maps satisfying a Markov property.  (Reference [1] should be consulted for a discussion of the poles encountered in the study of a Markovian map $f$).  Note that the studies of a.c.i.m. for maps of the interval, and of SRB measures for H\'enon-like maps, are typically based on perturbations of a map satisfying a Markov property (for the use of slightly more general Misiurewicz-type maps see [9], which also gives references to earlier work).  
\medskip
	The function $\Psi(\lambda)$ that we have encountered is related to the {\it susceptibility} $\omega\mapsto\Psi(e^{i\omega})$ giving the response of a system to a periodic perturbation.  The existence of a holomorphic extension of the susceptibility to the upper half complex plane is expected to follow from {\it causality} (causality says that cause preceeds effect, resulting in a {\it response function} $\kappa$ having support on the positive half real axis, and its Fourier transform $\hat\kappa$ extending holomorphically to the upper half complex plane).  A discussion of nonequilibrium statistical mechanics [7] shows that the expected support and holomorphy properties hold close to equilibrium, or if uniform hyperbolicity holds.  In the example discussed in this note, $\kappa$ has the right support property, but increases exponentially at infinity, and holomorphy in the upper half plane fails, corresponding the existence of a pole of $\Psi$ at $\lambda=1/\sqrt{f'(-1)}$.  This might be expressed by saying that $\rho$ is {\it not linearly stable}.  The physically interesting situation of {\it large systems} (thermodynamic limit) remains quite unclear at this point.  
\medskip\noindent
{\bf Acknowledgments.}  
For many discussions on the subject of this note, I am indebted to V. Baladi, M. Benedicks, G. Gallavotti, M. Viana, and L.-S. Young.  
\medskip\noindent
{\bf References.}

[1] V. Baladi, Y. Jiang, H.H. Rugh.  ``Dynamical determinants via dynamical conjugacies for postcritically finite polynomials.''  J. Statist. Phys. {\bf 108},973-993(2002).

[2] C. Bonatti, L. Diaz, and M. Viana.  {\it Dynamics beyond uniform hyperbolicity: a global geometric and probabilistic approach.} Springer, to appear.

[3] L. Chierchia and G. Gallavotti.  ``Smooth prime integrals for quasi-integrable Hamiltonian systems.''  Nuovo Cim. {\bf 67B},277-295(1982).

[4] J. P\"oschel.  ``Integrability of Hamiltonian systems on Cantor sets.''  Commun. in Pure and Applied Math. {\bf 35},653-696(1982).

[5] D. Ruelle.  ``Differentiation of SRB states.''  Commun. Math. Phys. {\bf 187},227-241(1997); ``Correction and complements.''  Commun. Math. Phys. {\bf 234},185-190(2003).  

[6] D. Ruelle.  ``Differentiation of SRB states for hyperbolic flows.''  In preparation.

[7] D. Ruelle.  ``Smooth dynamics and new theoretical ideas in nonequilibrium statistical mechanics.''  J. Statist. Phys. {\bf 95},393-468(1999).  

[8] D. Ruelle.  ``Application of hyperbolic dynamics to physics: some problems and conjectures.''  Bull. Amer. Math. Soc. (N.S.) {\bf 41},275-278(2004).

[9] Q. Wang and L.-S. Young.  ``Towards a theory of rank one attractors.''  Preprint

[10] H. Whitney.  ``Analytic expansions of differentiable functions defined in closed sets.''  Trans. Amer. Math. Soc. {\bf 36},63-89(1934).

[11] L.-S. Young.  ``What are SRB measures, and which dynamical systems have them?''  J. Statist. Phys. {\bf 108},733-754(2002).

\end